\documentclass[a4paper,12pt]{amsart}

\usepackage{amssymb,amscd,amsthm,epsfig,graphicx,amsmath}

\newcommand{\en}{\subset}

\newcommand{\Z}{\mbox{$\mathbb{Z}$}}

\newcommand{\R}{\mbox{$\mathbb{R}$}}
\newcommand{\T}{\mbox{$\mathbb{T}$}}

\newtheorem{teo}{Theorem}[section]
\newtheorem{cor}{Corollary}[section]
\newtheorem{lema}{Lemma}[section]
\newtheorem{prop}{Proposition}[section]

\newtheorem{defi}{Definition}[section]

\newcommand{\bi}{\begin{itemize}}
\newcommand{\ei}{\end{itemize}}

\theoremstyle{definition} \theoremstyle{remark}
\newtheorem{obs}{Remark}[section]

\newcommand{\D}{\mbox{$\mathbb{D}$}}

\newcommand{\eps}{\varepsilon}
\newcommand{\diam}{\hbox{diam}}

\newcommand{\interior}{{\rm int}}

\newcommand{\dem}{\vspace{.05in}{\sc\noindent Proof.} }
\newcommand{\demos}[1]{\vspace{.05in}{\sc\noindent Proof #1.} }
\newcommand{\lqqd}{\par\hfill {$\Box$}}

\newcommand{\umb}[1]{\hat #1}
\newcommand{\dimtop}{\hbox{dimtop}}
\newcommand{\cmax}{\mathcal{C}}
\newcommand{\dist}{\hbox{dist}}
\newcommand{\sing}{\mathcal{S}}

\author[A. Artigue]{Alfonso Artigue}
\address{Imerl, Facultad de Ingenier\'\i a, Universidad de la Rep\'ublica, Uruguay}
\email{aartigue@fing.edu.uy}
\author[J. Brum]{Joaquin Brum}
\address{Cmat, Facultad de Ciencias, Universidad de la Rep\'ublica, Uruguay}
\email{joaquin@cmat.edu.uy}
\author[R. Potrie]{Rafael Potrie}
\address{Imerl, Facultad de Ingenier\'\i a, Universidad de la Rep\'ublica, Uruguay}
\email{rpotrie@cmat.edu.uy}
\title{Local Product Structure for Expansive Homeomorphisms}

\begin{document}

\maketitle

\begin{abstract}
Let $f\colon M\to M$ be an expansive homeomorphism with dense
topologically hyperbolic periodic points, $M$ a closed manifold.
We prove that there is a local product structure in an open and dense
subset of $M$. Moreover, if some topologically hyperbolic periodic
point has codimension one, then this local product structure is
uniform. In particular, we conclude that the homeomorphism is
conjugated to a linear Anosov diffeomorphism of a torus.

\end{abstract}

\section{Introduction}

Let $M$ be a compact connected boundaryless manifold of dimension $n$ and $f: M
\to M$ an expansive homeomorphism, that is, there exists
$\alpha>0$ such that every two points have iterates which are
separated at least $\alpha$ from each other (the existence of
$\alpha$ is independent of the metric, furthermore, the notion can
be defined independently of the metric).

A paradigm of expansive homeomorphisms are Anosov diffeomorphisms.
Other class of expansive homeomorphisms are pseudoAnosov maps in
surfaces of genus $g>2$. They satisfy that $\Omega(f)=M$ and they
have dense topologically hyperbolic periodic points. In surfaces,
pseudoAnosov maps and linear Anosov homeomorphisms (that is,
conjugated to a linear Anosov diffeomorphism) describe completely
expansive dynamics as was proved in \cite{L2},\cite{H2} obtaining a global classification
of expansive homeomorphisms. For
this classification, the key step is to prove that in a reduced
neighborhood of every point there is a local product structure. To
do this, in \cite{L2} it is proved that every point in an
expansive homeomorphism has a uniformly big connected stable and
unstable set. On surfaces, in some way this is enough to find the
local product structure since proving that the connected sets
intersect is enough (using Invariance of Domain Theorem, see
\cite{Sp}) to find local product structure (these connected
sets contain arcs and so a map from $[0,1]^2$ to a neighborhood of
the point can be constructed). In higher dimensions, the existence
of connected stable and unstable sets is not enough to find a
local product structure, as shown in the example from
\cite{FR}.

A surprising result is the one of \cite{V2}, since it proves
that in dimension $3$ expansive homeomorphisms whose topologically
hyperbolic periodic points are dense, are conjugated to linear
Anosov diffeomorphisms in the torus $\T^3$. For doing this it is
also very important to find a local product structure in an open
and dense subset of the manifold (see \cite{V1}). Again, the
technique is to obtain intersections between stable and unstable
sets of topologically hyperbolic periodic points which are near
and use Invariance of Domain Theorem. This is not completely
direct since, a priori, the size of the stable and unstable sets
of the periodic points is not controlled, and must study
separation properties of these sets to ensure the intersection.
The hypothesis of having dense topologically hyperbolic periodic
points was weakened in \cite{V4} changing it for having
$\Omega(f)=M$ (a necessary condition as can be seen with the
example in \cite{FR}) and some smooth hypothesis ($f$ must
be a $C^{1+\theta}$ diffeomorphism) to use Pesin theory.

In this work we obtain local product structure in an open dense
subset of $M$ when topologically hyperbolic periodic points are
dense in $M$; in fact, we obtain local product structure in
neighborhoods of every periodic point. When the codimension of
topologically hyperbolic periodic points is arbitrary, this result
is optimal, since in the case of a product of two pseudoAnosov
maps the local product structure can not be defined in all the
manifold.

The somewhat strange aspect of the result from \cite{V2} is
that it proves that in dimension $3$ no singularities can appear,
not as in the surface case where pseudoAnosov maps are expansive
with dense topologically hyperbolic periodic points. However, this
result has a nice counterpart in the theory of Anosov
diffeomorphisms where it is known that codimension one Anosov
diffeomorphisms can only exist in torus and be conjugated to a
linear one (see \cite{F},\cite{N}).

Maybe this connection is not a priori obvious, but we give in this
work more evidence of it, proving that if the topologically
hyperbolic periodic points are dense in $M$ (with dimension higher
than $2$) and one of them has codimension one, then, the
homeomorphism is conjugated to a linear Anosov diffeomorphism of
$\T^n$. The reason why this does not work in dimension $2$ is that
we can disconnect an arc by removing from it one point and not a
disc of dimension $>1$. The proof in this case is based on proving
first that singularities are finite, and then discarding their
existence.

\subsection{Definitions and presentation of results}

In this section we define the concepts that we use in the course
of this paper and give precise statements of the results in it.

\begin{defi}
We say an homeomorphism $f\colon M\to M$ is \emph{expansive} if
$\alpha>0$ exists satisfying that if $x,y\in M$ are different
points, then, there exists $n\in \Z$ such that
$\dist(f^n(x),f^n(y))> \alpha$.
\end{defi}

\begin{defi} We say that a periodic point $p\in M$ of period $l$
is \emph{topologically hyperbolic} ($p\in Per_H$) if $f^l$ is
locally conjugated to the linear map $L\colon \R^r\times
\R^{n-r}\to \R^r\times \R^{n-r}$ given by $L(x,y)=(x/2,2y)$. In
this case we say that $p\in Per_H^r\subset Per_H$, we say that $r$
is the index of $p$.
\end{defi}

In our case $f$ is expansive, so, due to results in
\cite{L1} (Lemma 2.7) it is true that
$Per^0_H=Per^n_H=\emptyset$ since no stable points exist.

We denote as $H_k(A)$ ($H_{c}^k(A)$) the $k$ dimensional reduced
homology (cohomology with compact support) of $A$ with
coefficients in $\R$. As usual, we define the stable and unstable
sets of a point $x\in M$ as $W^s(x)=\{y\in
M:\dist(f^n(x),f^n(y))\to 0, n\to+\infty\}$ and $W^u(x)=\{y\in
M:\dist(f^n(x),f^n(y))\to 0, n\to-\infty\}$. The local stable and
unstable sets ($\eps$-local) are defined as follows
$W^s_\eps(x)=\{y\in M:\dist(f^n(x),f^n(y))\leq\eps,\forall n\geq 0
\}$ and $W^u_\eps(x)=\{y\in M:\dist(f^n(x),f^n(y))\leq\eps,\forall
n\leq 0 \}$. We denote as $cc_p(X)$ the connected component of
$X\subset M$ containing $p$.

We prove a separation property verified by the stable and unstable
set of a point $p\in Per_H^r$. The proof of this Proposition
follows the ideas in \cite{V1},\cite{V2}  and it is
developed in section \ref{seccionseparacion}. The property is the
following.

\begin{prop}\label{separacionestable}
  Let $f:M \to M$  be an expansive homeomorphism.
  Then, there exists $\eps>0$ such that for all
  $x\in M$,
  $p\in Per_H^k\cap B_\eps(x)$ and $V\subset B_\eps(x)$
  homeomorphic to
  $\R^n$ and containing $p$,
  we have
  $H_{n-k-1}(V\setminus S_p)\cong \R$ with $S_p=cc_p(V\cap
  W^s(p))$.
\end{prop}

An analogous result is verified for the unstable set.

\begin{obs}\label{observacionestablebola} If $f:M\to M$ is an expansive homeomorphism and $z \in
M$ then, for all $\eps>0$ exists $\delta>0$ such that if
$S_z=cc_z(W^s(z)\cap B_{\delta}(z))$ then $S_z \en W^s_{\eps}(z)$.
See \cite{L2}.
\end{obs}

\begin{defi}
  We say that $p\in M$ admits a \emph{local product structure} if
  there exists a map
  $h\colon \R^k\times \R^{n-k}\to M$
   which is a homeomorphism over its image
  ($p \in Im(h)$) and if there exists $\eps>0$
  such that for all $(x,y)\in \R^k\times \R^{n-k}$ it is verified
  that $h(\{x\}\times \R^{n-k})= W^s_\eps (h(x,y))\cap Im(h)$ and
  $h(\R^k\times\{y\})= W^u_\eps (h(x,y))\cap Im(h)$. We say that the local product structure is
  a \emph{uniform local product structure} if in addition to the previous conditions,
  there exists $r>0$ such that for all $x \in M$ the points in
  $B_{r}(x)$ admit a local product structure.
\end{defi}

We remark that the points admitting a local product structure are
an open set. We call the points which do not admit a local product
structure singularities.

\begin{teo}\label{teo1}
  Let $f\colon M\to M$ be an expansive homeomorphism such that
  $\overline{Per_H}=M$. Then, every point in $Per_H$ admits a local
  product structure. In particular, the set of points with a local
  product structure is open and dense in $M$.
\end{teo}

Once this is obtained, in \cite{V2} singularities are studied,
discarding their existence by studying the way in which the
product structure is glued together in the singularity and proving
that this can not happen by discarding  the possible dimensions in
which that gluing may happen one by one . As was already
explained, with the product between the Anosov and the
pseudo-Ansov we see that this can not be done in dimension larger
than $3$, unless we add the hypothesis of having
$Per_H^{n-1}=Per_H$. This will be studied in section
\ref{seccionsingularidades}.

It is worth observing that the fact of having a local product
structure in an open and dense subset does not imply, a priori,
that the index of the topologically hyperbolic periodic points
should be constant in all the manifold. We shall prove this is
true, under the hypothesis of Theorem \ref{teo1}, for dimensions 3
and 4. For doing that, in section \ref{seccionlemas} several
properties of the points in $\overline{Per_H^{n-1}}$ are studied.
 The following sharper result is obtained.

\begin{teo}\label{teo2}
Let $f:M\to M$ be an expansive homeomorphism verifying
$\overline{Per_H} = M$. Then, $Per_H^{n-1}=Per_H$ or
$Per_H^{n-1}=\emptyset$. Analogously for $Per_H^1$.
\end{teo}

\begin{cor} Let $f:M \to M$  be an expansive homeomorphism of a manifold
of dimension $3$ or $4$ with $\overline{Per_H(f)} =M$. Then, every
topologically hyperbolic point has the same index.
\end{cor}

\dem{} In dimension $3$ we have $Per_H=Per_H^1\cup Per_H^2$ (see
\cite{L1}, Lemma 2.7, no stable points can exist); the
Theorem \ref{teo2} concludes the proof. In dimension $4$, we have
$Per_H=Per_H^1\cup Per_H^2\cup Per_H^3$ and since $Per_H^1\cup
Per_H^3=\emptyset$ implies $Per_H=Per_H^2$ the proof finishes by
using the Theorem \ref{teo2}. \lqqd

Finally, in section \ref{seccionsingularidades} we study the
singularities in the case of having one topologically hyperbolic
point of index $n-1$ (or $1$), discarding their existence and
concluding that there is a uniform local product structure in all
the manifold.

\begin{defi} Let $f:M\to M$ be a homeomorphism, we say it verifies
the \emph{pseudo orbit tracing property} if for all $K>0$
exists $\alpha>0$ such that if $\{x_n\}_{n\in \Z}$ verifies $\dist(x_n,
f(x_{n-1}))< \alpha$ (i.e. it is an $\alpha-$pseudo-orbit) then
there exists $x \in M$ such that  $\dist(f^n(x), x_n) < K$ for all
$n \in \Z$ (i.e. $x$ $K-$shadows the pseudo orbit).
\end{defi}

\begin{teo}\label{teo3}Let $f:M\to M$ be an expansive homeomorphism
verifying $\overline{Per_H} = M$ and $Per^{n-1}_H \neq \emptyset$
or $Per^{1}_H \neq \emptyset$ . Then, there is a uniform local
product structure in all the manifold. In particular, the pseudo
orbit tracing property is verified.
\end{teo}

In dimension $3$, in \cite{V3} the uniform local
product structure is used for proving that $M=\T^3$ and concluding
that $f$ must be conjugated to a linear Anosov diffeomorphism. In
higher dimensions, as far as we know, there are no published
results which ensure that a manifold with uniform local product
structure of codimension one is a torus. However, our results give
a codimension one foliation transversal to a dimension one
foliation. It is known from the work of Franks that if the
foliations are differentiable this implies that the manifold is a
torus. This is also the case without the differentiability
assumption. The proof is a straightforward adaptation of the work
in \cite{V3} and \cite{F}. However, we shall
sketch how to adapt the proof for the sake of completeness. We
then have the following Corollary, which is the main result of
this paper.

\begin{cor}Let $f:M^n \to M^n$ ($n \geq 3$) be an expansive homeomorphism verifying
$\overline{Per_H} = M$. Suppose $Per^{n-1}_H \neq \emptyset$ or
$Per^{1}_H \neq \emptyset$. Then, $M= \T^n$ and $f$ is conjugated
to a linear Anosov diffeomorphism.
\end{cor}

\dem{} It is consequence of the Theorem \ref{teo3} and a result of Hiraide (\cite{H1})
which ensures that an expansive homeomorphism in $\T^n$ with the
pseudo orbit tracing property is conjugated to a linear Anosov
diffeomorphism. The proof that $M=\T^n$ is sketched at the
Appendix of this work.

\lqqd

\emph{Acknowledgements:} We profoundly thank  Jos\'e Vieitez for his
patience and suggestions during the preparation of this work. We
want to thank him and Jorge Lewowicz for teaching us all we know
about this topics. Also, we thank Alberto Verjovsky for his
suggestions and comments on topology (especially for communicating
to us the proof sketched in the Appendix) . Finally, we would like
to thank the members of the dynamical systems group in Montevideo
for useful conversations and advice.


\section{Separation properties}\label{seccionseparacion}

In this section, with the help of the ideas in \cite{V1}, we
prove the Proposition \ref{separacionestable}.

The following Lemma is a general homological property of euclidean
spaces.

\begin{lema}\label{separacionabstracta} Let $B$ be a set homeomorphic to $\R^n$ and $F \en B$
a closed connected set homeomorphic to an open set of $\R^k$.
Then, $H_{n-k-1}(B\setminus F) \cong \R$.
\end{lema}

\dem{} Let $U=B \setminus F$. We then have the following long
exact sequence of homology:
    \[
        \ldots \to H_l(B) \to H_l(B,U) \to H_{l-1}(U) \to H_{l-1}(B)
        \to \ldots
    \]
    We know $H_l(B)=0$ (recall we work with reduced homology), so, we have $H_l(B,U) \cong H_{l-1}(U)$.

    In particular, it is true that $H_{n-k}(B,U) \cong
    H_{n-k-1}(U)$. Using the duality Theorem of Alexander-Pontryagin
    we also deduce that $H_{n-k}(B,U)\cong H_c^{k}(F)$ (see \cite{D},\cite{Sp}). Applying the same Theorem,
    now to $(F,\emptyset)$, we can conclude that $H_{0}(F, \emptyset) \cong
    H_{c}^k(F)$. Therefore, we can deduce (using that $H_{0}(F,
    \emptyset)\cong \R$ since $F$ is connected) that $H_{n-k-1}(U) \cong \R$ as we desired.

\lqqd

\begin{lema}\label{compconestable} Let $f:M \to M$  be an expansive homeomorphism.
Then, there exists $\eps>0$ such that for all $p \in Per^r_H$
exists $\phi: \overline{D^r} \to W^s(p)$ a surjective
homeomorphism over its image satisfying that: $\phi(0)=p$ and that
for all continuous curve $y: [0,1] \to \overline{D^r}$  such that
$y(0)=0$ and $y(1) \in
\partial D^r$ there exists $s \in (0,1]$ such that $\phi \circ y (s)
\notin B_{\eps}(p)$.
\end{lema}

\dem{} Expansivity ensures the existence of $\eps>0$ such that for
every connected set $C$ with diameter smaller than $\eps$
satisfying that the diameter of $f^n(C)$ is bigger than the
constant $\alpha$ of expansivity for some $n \leq 0$, then, the diameter of
$f^m(C)$ is bigger than $\eps$ for all $m<n$.

If this affirmation were false there would exist connected sets
$C_n$ with diameter smaller than $1/n$ and numbers $k_n>0$ and
$l_n>k_n$ verifying that the diameter of $f^{-k_n}(C_n)$ is bigger
than the expansivity constant and the diameter of  $f^{-l_n}(C_n)$
smaller than $1/n$. Using the uniform continuity of $f$ we obtain
that $k_n \to +\infty$ and $l_n-k_n \to +\infty$. Connectedness of
$C_n$ and its iterates allows us to find points $x_n$ and $y_n$ in
$f^{-m_n}(C_n)$ (with $0 \leq m_n \leq l_n$ and $l_n -m_n \to
\infty$) such that $\alpha/2 \leq d(x_n,y_n) < \alpha$ and $d(f^i
(x_n), f^i(y_n))< \alpha$ for all $- l_n + m_n \leq i \leq m_n$.
Considering limit points of the sequences $x_n$ and $y_n$ we
contradict the expansivity of $f$.

Without loss of generality we can suppose that $p$ is a fixed
point and we can consider the conjugation $h : \overline{D^r} \to
W^s(p)$ between $f$ and the linear hyperbolic map. Also, we know
that there exists $N<0$ such that for all  $x\in h(\partial
D^r)\subset W^s(p)$ exists $n\in [N,0]$ satisfying $f^n(x)\notin
B(p,\alpha)$ (if not, we can find points in $h(\partial D^r)$
which stay in $B(p,\alpha)$ for an arbitrarily large quantity of
iterates of $f$, taking limit points of that sequences we
contradict expansivity). We then define $\phi: \overline{D^r} \to
W^s(p)$ by $ \phi(x)=f^N \circ h(x)$. Then, for every $y$
connecting $p$ with $\phi(\partial D^r)$ we have that $y([0,1])$
is a connected set of diameter bigger than $\eps$. For this $\eps$
the lemma works. \lqqd

\demos{of Proposition \ref{separacionestable}} After what we have
already proved, to conclude the proof, it is enough to prove that
if we have a homeomorphism over its image $\phi:\overline{D^k} \to
\R^n$ such that $\phi(0)=0$ and such that for every curve $y:[0,1]
\to \overline{D^k}$ verifying $y(0)=0$ and $y(1) \in
\partial D^k$ satisfies that  $\phi\circ y ([0,1])$ is not contained in $B_\eps(0)$,
so, considering $X$, the connected component of $0$ in
$\phi^{-1}(B_\eps(0))$ we have $H_{n-k-1}(B_\eps(0)\setminus
\phi(X))= \R$.

In order to do this, let $F = \phi(X)$ and $B=B_\eps(0)$. Since
$B$ is open, we have that $\phi^{-1}(B)$ is an open set of
$\overline{D^k}$. Since $\overline{D^k}$ is locally arcconnected,
$X$, being a connected component of an open set is open in
$\overline{D^k}$ and locally arcconnected. This implies that it is
arcconnected.

We have that $X \cap
\partial D^k = \emptyset$ since in the other case a curve joining
$0$ with $\partial D^k$ whose image by $\phi$ would be included in
$B$ would exist. Then, $F$ is homeomorphic to an open set of
$\R^k$. Since $X$ is a connected component, $X$ is closed in
$\phi^{-1}(B)$ so $F$ is closed in $B$. Lemma
\ref{separacionabstracta} implies the thesis. \lqqd


\section{Local product structure}

    The construction of a local product structure is strongly
    based on proving that stable and unstable sets of the periodic
    points intersect. This allows us to define a map between
    $W^s_{\eps}(p) \times W^u_{\eps}(p)$ and a neighborhood of
    $p$ which is a homeomorphism by the invariance of domain theorem and has
    the desired properties. In this section we prove that this
    intersection occurs for periodic points close to a given one.

  Let $\Delta^{m}=\{ (x_1, \ldots x_{m+1}) \in \R^{m+1} \ :
  \ x_i \geq 0 \ , \ x_1 + \ldots x_{m+1} =1 \}$ the canonical simplex of dimension $m$. We denote
  $\gamma = \sum_i a_i \sigma_i$
  to a $m-$chain, where $\sigma_i \colon
  \Delta^{m}\to V$ ($a_i\in\R$). In the course of this section, $\gamma$
  denotes the chain and the union of the images of $\sigma_i$
  indifferently.

\begin{lema}\label{prodloc1}

  For all $x\in M$, there exists $\eps>0$ such that if $V\subset
  B_\eps(x)$ is homeomorphic to $\R^n$ and $p\in V\cap Per_H^l$
  then there exists a cycle $\gamma \en U_p$ which is non trivial in the
  $n-l-1$ dimensional homology of $V\setminus S_p$
   (where $S_p=cc_p(V\cap
  W^s(p))$ and $U_p=cc_p(V\cap
  W^u(p))$). Furthermore,  given $K$ compact in $V$ we can choose
  $\gamma$ so that $\gamma \en V \setminus K$.
\end{lema}

\begin{dem}
  Because of Proposition \ref{separacionestable} we know that
  $\eps_0>0$ exists verifying that $H_{n-l-1}(V\setminus S_p)\neq 0$.
  Let $\gamma$ be a cycle such that its $n-l-1$ dimensional homology class
  $[\gamma]$ is non trivial. Since $H_{n-l-1}(V)=0$ we can suppose $\gamma=\partial \eta$
  where $\eta$ is a $n-l$ dimensional chain in V.

  Say $\eta=\sum_{i=1}^{j} a_i\sigma_i$ with $\sigma_i\colon
  \Delta^{n-l}\to V$ ($a_i\in\R$).

  Besides, we can suppose that $\sigma_i$ and $\partial \sigma_i$ are
  topologically transversal to $S_p$ so that the set of points of intersection between
  every $\sigma_i$ and $S_p$ is finite and such that $\partial \sigma_i\cap
  S_p=\emptyset$. Given $\eps_1>0$, using barycentric subdivision (see \cite{Sp}), we can also suppose
  $\diam(\sigma_i)<\eps_1$. We observe that if
  $\sigma_i\cap S_p=\emptyset$ then $\partial \sigma_i$ is trivial in $H_{n-l-1}(V\setminus S_p)$. So, by choosing $\eps_1$
  small enough we can suppose that each $\sigma_i$ intersects
  $S_p$ in $y_i$ only for $i=1,...,j$.

    Let $h\colon U\subset \R^n\to M$ the local conjugation with
    the hyperbolic map, in a neighborhood of $p$. Intersecting with $V$ we have
    that
    $h(U)\subset V$ and by iteration of $f$ we
    can suposse that is a neighborhood of $S_p$.

    We can think $U\subset V\subset
    \R^n$ (with the identification given by $h$)
    , $S_p\subset \R^l\times\{p_2\}$ and $U_p\subset \{p_1\}\times
    \R^{n-l}$ where $p=(p_1,p_2)$.

  We can choose $\eps_1$ smaller so that $B_{\eps_1}(y_i)\subset U$

  Since $y_i \in \sigma_i$ and $\diam(\sigma_i)< \eps_1$ we have $\sigma_i\subset B_{\eps_1}(y_i)$.
  Let $h^i_t\colon \R^l\times\R^{n-l}\to
  \R^l\times\R^{n-l}$ continuous given by $h^i_t(a+y_i^1,b)=(ta +y_i^1,b)$ with $t\in [0,1]$
  where $y_i=(y_i^1,y_i^2)$.

  Then, for $t\in[0,1]$, $h^i_t\circ \partial\sigma_i$ does not intersect $S_p$ and is contained in $V$.
  Also, we have $h^i_1\circ \partial\sigma_i = \partial\sigma_i$ and $h^i_0\circ \partial\sigma_i
  \subset \{y_i^1\}\times \R^{n-l}$.
  Since $h_0\circ \partial\sigma_i$
  is homotopic to
  $\partial\sigma_i$ we have they are both homologous in
  $V\setminus S_p$.

  For every $i=1,...,j$ let
  $\beta_i\colon [0,1]\to S_p$ be a continuous curve such that $\beta(0)=y_i$ and $\beta(1)= p$.
  If we choose a smaller $\eps_1$ again, we have
  $B_{\eps_1}(\beta_i)\subset U$ for all $i=1,...,j$.

  Now, we consider
  $g^i_t\colon\R^n\to\R^n$, another homotopy, given by
  $g^i_t(z)=z+\beta_i(t)-y_i$. It verifies that $g^i_t(h_0\circ \partial\sigma_i)$
  does not intersect $S_p$ for all $t\in[0,1]$, $g^i_0=id_{\R^n}$
  and  $g^i_1(y_i)=p$ so $\sum_{i=1}^{j}a_i g^i_1\circ
  h^i_0\circ\partial\sigma_i \subset U_p$, and since $g_t^i$ is a homotopy, it is homologous to
  $\gamma=\sum_{i=1}^{j}a_i\partial\sigma_i$ which is non trivial in the homology of
  $V\setminus S_p$. We call $\gamma$ to $\sum_{i=1}^{j}a_i g^i_1\circ
  h^i_0\circ\partial\sigma_i$.

  \medskip

  To see that there is a cycle homologous to $\gamma$ outside of every compact set in $V$,
  we will use the map of the Lemma \ref{compconestable}
  $\phi\colon \overline{D^{n-l}}\subset \R^{n-l}\to
  M$ which verifies that $U_p=\phi (X)$ where
  $X=cc_0(\phi^{-1}(V))$. Consider a subdivision of $\R^{n-l}$ in
  simplexes of  dimension $n-l$ and diameter smaller than $\rho$.
  Let us say  $\R^{n-l}=\bigcup_{i=1}^\infty \theta_i$ and that
  $0\in\R^{n-l}$ is in the interior of $\theta_0$.

    If we consider a neighborhood $B\subset V$ of $p$ with linear
    structure as before, we know that $H_{n-l-1}(B\setminus
    S_p)\cong \R$. So, we have that there exists a non zero $a\in
    \R$ such that
    $\gamma=a\partial(\phi\circ\theta_0)$ in $H_{n-l-1}(B\setminus
    S_p)$ and in particular also in $H_{n-l-1}(V\setminus
    S_p)$.
  Let $\eta_1=\theta_0-\sum_{\theta_i\subset X}\theta_i$.

  We observe that $\partial(\phi\circ\eta_1)$ is a
  trivial cycle in $V\setminus S_p$.
  So, $a^{-1} \gamma$ is homologous to
  $\gamma'=\partial\phi\circ(\sum_{\theta_i \en X}\theta_i)=
  \sum_{\theta_i\en X}\phi\circ\partial\theta_i$.

  To conclude the proof is enough to observe that we can suppose $\sum_{\theta_i\en X}
  \partial\theta_i \en B_{\rho}(\partial X)$ and use the fact
  that $\phi$ is uniformly continuous. This is true because every boundary in
  $B_{\rho}(\partial X)$ is cancelled for being trivial in homology and we can take
  $\theta_i$ to have arbitrarily small diameter. Given a compact
  set in $V$, considering an adequate $\rho$ we conclude the proof.
\end{dem}\lqqd

\begin{cor}\label{sepcoduno}
  With the same hypothesis that the previous Lemma, if $p\in Per_H^{n-1}$
  then $S_p$ separates $V$ in two connected components  $V_1$ and $V_2$.
  Also, $p$ separates $U_p$ in two connected components $U_1$ and
  $U_2$ such that $U_1\subset V_1$ and $U_2\subset V_2$.
\end{cor}

\dem
  Due to the fact that we are working with reduced homology, the previous Lemma
  implies that
  $V\setminus S_p$ has two connected components $V_1$ and $V_2$. Moreover, $U_p$
  is homeomorphic to $\R$, so $U_p\setminus \{p\}$ has two connected components $U_1$ and $U_2$.
  Let us suppose that $U_1,U_2\subset V_1$. Since $V_1$ is connected, we have that every
   $\gamma\subset U_1\cup U_2$ would be trivial in the homology of $V\setminus S_p$,
  contradicting the previous Lemma.
\lqqd

We will repeatedly make use of the following Lemma concerning the
semicontinuous variation of stable and unstable sets (see
\cite{L2}).

\begin{lema}\label{semicontinuidad} Let $f: M \to M$ be any homeomorphism.
Then, given $\eps, \gamma >0$ and $x \in M$, there exists $\delta>0$
verifying that if $\dist(x,y)< \delta$ then, $W^{s}_{\eps}(y) \in
B_{\gamma}(W^s_{\eps}(x))$. \end{lema}

\dem{} Suppose by contradiction that there exists $\gamma,\eps>0$ and $x_n \to x$ such
that $y_n \in W^s_{\eps}(x_n) \cap B_{\gamma}(W^{s}_{\eps}(x))^c$
exist. If we consider $z$ a limit point of $y_n$ we have
$$\dist(f^{k}(z), f^{k}(x))= \lim_{n \to +\infty}
\dist(f^{k}(y_n), f^k(x_n)) \leq \eps$$ with $z\neq x$ and $k\geq 0$. Thus, $z\in W^s_\eps(x)$, but this is a contradiction since $z \notin B_{\gamma}(W^s_\eps(x))$. \lqqd

Another result we will repeatedly make use of refers to the
distance between local stable and unstable sets of the points (see
\cite{V1} also). We think of it as ensuring ``big angles''
between the local stable and unstable sets.

\begin{lema}\label{angulos}Let $f:M \to M$ be an expansive homeomorphism with expansivity constant $\alpha>0$.
Given $V \en U$ neighborhoods of $x$ and $\rho$ small enough, there exist a neighborhood $W \en
V$ of $x$ such that if $y, z \in W$ we have $\dist(S_y \cap (U
\setminus V), U_z \cap (U \setminus V)) > \rho$ (where $S_y= cc_y
(W^s(y) \cap U$ and $U_z =cc_z (W^u(z) \cap U$).
\end{lema}

\dem{} For $0<\eps<\alpha$ let us consider $\delta>0$ given by Remark \ref{observacionestablebola}. Then, we can see that for given neighborhoods $V\en U \en B_\delta(x)$ of $x$, there are $\rho >0$ and and $W \en V$ ($x\in W$) such that if $y,z \in W$, then $$dist (S_y \cap (U\backslash V), U_z \cap (U \backslash V)) >\rho$$

 Otherwise, there would be points $y_n$ and $z_n$ converging
to $x$ and such that $\dist(S_{y_n} \cap (U \setminus V), U_{z_n}
\cap (U \setminus V) )< 1/n$. Taking a limit point of $a_n \in
S_{y_n} \cap (U \setminus V)$ (choosen to verify
$\dist(a_n,U_{z_n} \cap (U \setminus V))< 1/n$) we find a point
$\overline{x}\neq x$ such that $\overline{x} \in S_x \cap S_y \cap \overline{(U \backslash V)}$. Thus, by Remark \ref{observacionestablebola} $$\dist(f^k(x), f^k(\overline{x})) \leq
\eps < \alpha$$ $\forall k \in \Z$ so, expansivity implies $x=\overline{x}$ which is a contradiction.

\lqqd


In the following Proposition we prove that the index of
topologically hyperbolic periodic points is locally constant and
that if two of them are close enough then their local stable and
unstable sets intersect. As was already mentioned, this is the key
step for obtaining the local product structure.

\begin{prop}\label{nosequeponer}
  Let $f\colon M\to M$ be an expansive homeomorphism. Then
  \begin{enumerate}
    \item for all  $k=1,...,n-1$, $Per_H^k$ is open in $Per_H$ and
    \item for all  $p\in Per_H$ there exists open neighborhoods of $p$, $V_1$ and $V_2$
    such that for all $q\in Per_H\cap V_1$ we have $S_q\cap U_p\neq\emptyset$ and $U_q\cap
    S_p\neq\emptyset$, where $S_x=cc_x(W^s(x)\cap V_2)$, $U_x=cc_x(W^u(x)\cap
    V_2)$.
  \end{enumerate}
\end{prop}

\begin{dem}
  Let $p\in Per_H^k$, $\eps>0$ from Lemma \ref{prodloc1} applied to $p$ and $h\colon
  B_\rho(0)\subset \R^n\to h(B_\rho (0))\subset B_{\eps}(p)$ the local conjugacy, $h(0)=p$,
  between $f$ and $L\colon \R^k\times \R^{n-k}\to\R^k\times
  \R^{n-k}$ given by $L(x,y)=(x/2,2y)$, considering in $\R^n=\R^k\times
  \R^{n-k}$ the metric
  $d((x,y),(u,v))=\max\{\|x-u\|,\|y-v\|\}$. Fix $\rho_1\in
  (0,\rho)$  and let $V_2=h(B_{\rho_1}(0))$.
  For $q\in h(B_{\rho_1}(0))$ we denote $S'_q=h^{-1}(S_{q})$
  and $U'_q=h^{-1}(U_{q})$.

  Let $\rho_2$ and $\rho_3$ given by Lemmas  \ref{semicontinuidad} and
 \ref{angulos} such that if $\dist(h^{-1}(q),0)<\rho_3$ then
  $U'_q\cap B_{\rho_2}(\overline{S'_p}\cap\partial
    B_{\rho_1}(0))=\emptyset$ and
  $S'_q\subset B_{\rho_2}(S'_p)$.
  Let $V_1= h(B_{\rho_3}(0))$. Observe that we can use Lemma
  \ref{semicontinuidad} to $S_p$ and $U_p$ because of the choice
  of $\rho_1$.

  By applying Lemma
  \ref{prodloc1}  we know that if $q\in V_1\cap Per^m_H$ then there exists
  $h\circ\gamma \subset S_q$ a non trivial cycle of the $m-1$ dimensional homology of
  $V_2\setminus U_q$. Because of Lemma \ref{prodloc1} as well, we
  can suppose that
  $\gamma\subset B_{\rho_2}(\overline{S'_p}\cap \partial
  B_{\rho_1}(0))$.

  Let $\pi_t\colon \R^k\times \R^{n-k}\to \R^k\times
  \R^{n-k}$ given by $\pi_t(x,y)=(x,ty)$ for $t\in [0,1]$. It is easy to see that
  $\pi_t(B_{\rho_2}(\overline{S'_p}\cap\partial B_{\rho_1}(0)))\subset
  B_{\rho_2}(\overline{S'_p}\cap\partial
  B_{\rho_1}(0))$ for all $t\in[0,1]$. Then, $\pi_t\circ\gamma$ is
  a homotopy between
  $\gamma$ and $\pi_0\circ \gamma\subset S'_p$ contained in $B_{\rho_1}(0)
  \setminus U'_q$, so they are homologous in $B_{\rho_1}(0)\setminus U'_q$. To conclude:
  \begin{enumerate}
    \item If $Per_H^k$ is not open in $Per_H$ we can suppose there
    exists $q\in V_1\cap Per_H^m$ with $m<k$. Then, $\gamma$
    has dimension $m-1<k-1$ but $m-1$ dimensional homology of
    $S'_p\setminus U'_q$ is trivial (remember $S'_p$ is a disk) which is absurd.

    \item If $m=k$  a cycle $\eta\subset S'_p$ such that $\partial \eta = \gamma$ exists.
    Since $h\circ\gamma$ is non trivial in
    $V_2\setminus U_q$ we conclude that $S_p\cap U_q\neq\emptyset$.
  \end{enumerate}
\end{dem}\lqqd

\demos{of Theorem \ref{teo1}}

We are going to construct a local product structure in a
neighborhood of every $p \in Per_H$. We consider the notation of
the statement 
of Proposition \ref{nosequeponer}.

  Let $\pi_s\colon \overline{V_1} \to S_p$ be defined in the points $q\in
  Per_H$ as $\pi_s(q)=U_q\cap S_p$. This map is well defined in a dense subset of $V_1$ because
  of Proposition \ref{nosequeponer}.
  Let $x\in\overline{V_1}$ and $q_n\to x$, $q_n\in Per_H$ with $\pi(q_n)\to
  y$. Observe that $y \in W_{\varepsilon}^u(x) \cap S_p$ and
  expansivity imply that the intersection point is unique. This allows us to extend
  $\pi_s$ to $\overline{V_1}$. The same reason ensures this extension is continuous.
  Also we have $\pi_s(x)\in U_x\cap
  S_p$ and because of expansivity $\pi_s(x)= U_x\cap
  S_p$ for all $x\in V_1$. Expansivity also implies that $\pi_s|S_x$
  is injective.

  If $q\in Per_H$ then the Invariance of Domain Theorem  (see \cite{Sp}) implies that
  $\pi_s|S_q$ is open and a homeomorphism over its image.
  Observe that $\pi_s(r)\in \pi_s(S_q)$ with $r,q\in Per_H$
  implies $U_r\cap S_q\neq \emptyset$.
  Let $W\subset S_p$, $p\in W$, $W$
  homeomorphic to the disk $\overline D^k$ and $W$ relative
  neighborhood of
  $p$ in $S_p$.\\

  We affirm there exists $V_3$ neighborhood of $p  \ $ such that for
  all $q \in V_3 \cap Per_H$, $W\subset\pi_s(S_q)$. Otherwise, $q_n\to p$ would
   exists, such that $q_n\in Per_H$
  and $W\nsubseteq\pi_s(S_{q_n})$.
  Since $W$ is connected and $\pi_s|S_{q_n}$ open,
   $y_n\in \partial \pi_s(S_{q_n})\cap W$ must exist
  (the frontier is relative to $S_p$). So there must exist $x_n\in\partial V_1\cap S_{q_n}$
  such that $\pi_s(x_n)=y_n$.
  We can suppose $x_n \to x$ and $y_n \to y$ points of $S_{q_n}\cap\partial V_1$
  and $W$ respectively, the first
  due to semicontinuity of local stable sets (Lemma \ref{semicontinuidad})
  and the second because $W$ is compact.
  From the construction of $\pi_s$ we deduce that $x$ and $y$ are over the same
  local stable and unstable set contradicting expansivity
  (observe that $\dist$($W$,$S_p\cap\partial V_1)>0$ so $x \neq y$).\\

  Let $V_4=\pi_s^{-1}(W)\cap V_3$, we have that for every $q,r \in Per_H\cap
  V_4$, it is true that $S_q\cap U_r\neq \emptyset$ and $S_r\cap U_q\neq
  \emptyset$ from its construction.

  Let $A_s\subset S_p\cap V_4$
  and  $B_u\subset U_p\cap V_4$  be relative neighborhoods of $p$, both homeomorphic to disks.
  Now, let $x\in A_s$ and $y\in B_u$, then, by taking limit points of the intersection
   of local stable and unstable sets of periodic points converging to $x$ and $y$ respectively,
  semicontinuity of local stable and unstable sets (Lemma \ref{semicontinuidad})
  and expansivity easily imply that $U_x\cap S_y$ is a unique point. Let $h\colon
  A_s\times B_u\to V_1$ given by $h(x,y)=U_x\cap S_y$. It is
  continuous and injective. Using the Invariance of Domain Theorem again we conclude
  that it is open. This concludes the proof of the existence of a local product
  structure in an open and dense set.
\lqqd

\begin{obs} Although this does not ensure the dimension of
the decomposition in the local product structure to be constant,
it is an inmediate consequence of the obtained results if the
hypothesis of $f$ being transitive is added. We prove in section
\ref{secciondimconst} that the splitting is constant when
$Per_H^{n-1} \neq \emptyset$ or $Per_H^{1} \neq \emptyset$.
\end{obs}


\section{Codimension one case}

\subsection{Periodic point ordering and its properties}\label{seccionlemas}

We shall study the structure of $Per_H^{n-1}$ in a neighborhood of
a singularity $x \in M$ defining a partial order in $Per_H^{n-1}$.
We consider $B_{\nu}(x)$ so that Proposition
\ref{separacionestable} holds. Let

$$S_p = cc_p(W^s(p) \cap B_{\nu}(x)),$$
$$U_p = cc_p(W^u(p) \cap B_{\nu}(x)).$$

Also, we shall suppose that, because of Remark
\ref{observacionestablebola},  $S_p\subset W^s_\eps(p)$ and
$U_p\subset W^u_\eps(p)$ for some $\eps>0$. For every $p \in
Per_H^{n-1} \cap B_{\nu}(x)$ we define $\umb{p}= B_{\nu}\setminus
cc_x(B_{\nu}(x)\setminus S_p)$. 

Given $\delta>0$ we define the following order relation in
$X_\delta= Per_H^{n-1} \cap B_\delta (x)$. If $p,q \in X_\delta$
we say that $p \leq q$ if $\umb{p} \en \umb{q}$. Clearly this is a
partial order which depends on the singularity $x \in M$, $\nu>0$
from Proposition \ref{separacionestable} and $\delta \in (0,\nu)$.
We call chain to every totally ordered subset of the relation.

\begin{obs}\label{ordenarbol}
 Since stable sets of different periodic points have empty intersection, we have that if $\umb{p} \cap \umb{q}\neq \emptyset$ then the points $p$ and $q$ must be related by the ordering. So, if $p\leq q$ and $p\leq r$ then $\hat p\subset
\hat q \cap \hat r$, $q$ and $r$ must be related. This implies that if
$p\leq q$ and $\mathcal{C}$ is a maximal chain containing $p$, then $q\in \mathcal{C}$.
\end{obs}

This order can be well understood in the case of surfaces where, for the
pseudo Anosov maps, singularities have more than 2 maximal chains.

\begin{lema}\label{clasesdisjuntas}
  Given a singularity  $x\in M$ and $\nu >0$
  there exists $\delta>0$ such that
  there are finitely many maximal chains in $X_\delta$. These are pairwise disjoint
  and every one of them accumulates
  in $x$.
\end{lema}

\dem{} Let us suppose there were infinitely many maximal chains different
from each other. We shall prove this implies the existence of
arbitrarily large sets of points which are not pairwise related by the order
relation. We prove this using induction.

Let $p_1, ... , p_l \in X_\delta$ be pairwise not related. Let $\mathcal{C}_i$ be maximal chains such that $p_i\in \mathcal{C}_i$ and take $\mathcal{C}\neq \mathcal{C}_i$ another maximal chain. Since two points in the same maximal chain are related, at most one of the $p_i$'s can belong to $\mathcal{C}$.

If $p_i\notin \mathcal{C}$ for all $i=1,\dots,l$ then, we can choose $p_{l+1}\in \mathcal{C}\backslash \left(\bigcup_i \mathcal{C}_i\right)$ and it will not be related to any of the $p_i$ by Remark \ref{ordenarbol}.

If $p_i\in \mathcal{C}$ for some $1\leq i \leq l$ then, we can take $p'_i \in \mathcal{C}_i\backslash \mathcal{C}$ and $p_{l+1} \in \mathcal{C}\backslash \mathcal{C}_i$ not related. So, the points in $\{p_1, \ldots, p'_i, \ldots p_l, p_{l+1}\}$ will be pairwise not related again by Remark \ref{ordenarbol}.

This leads us to a contradiction since Lemma \ref{angulos} implies
the existence of $\delta>0$ and $\nu'\in (0,\nu)$ such that if
$p,q \in X_\delta$ then
$$
  \dist(S_p \cap \partial B_{\nu'}(x), U_q \cap\partial B_{\nu'}(x)) >
  \rho
$$

Given $p_i \in X_\delta$, Lemma \ref{prodloc1} ensures the
existence of $q_i\in\umb{p}\cap\partial B_{\nu'}(x)\cap U_{p_i}$
so that $\dist(q_i,q_j)>\rho$ if $i\neq j$. So, there exists a
bound on the number of pairwise not related points since $\partial
B_{\nu'}(x)$ is compact.

Once we know there are finitely many maximal chains, we know that the
ones that do not accumulate in $x$ are at a positive distance of
$x$, so if we choose $\delta$ to be smaller, we obtain that every
maximal chain in $B_{\delta}(x)$ accumulates in $x$.

Let $C$ and $C'$ be two maximal chains and $q \in C\cap C'$. If we
choose $\delta$ smaller in such a way that $\umb{q}$ be disjoint
with $B_\delta(x)$, we reduce the number of maximal chains in
$B_\delta (x)$. So, we can suppose that the maximal chains are
pairwise disjoint. \lqqd

We call $[p]$ to the maximal chain of $p$ in $X_\delta$ given by
the previous Lemma. Now, we define
$$
  S_{[p]}= \overline{\bigcup_{q \in [p]} \umb{q}} \en
  \overline{B_{\nu}(x)}
$$
\noindent where $p \in X_\delta$.

We remark that we can choose $\nu$ such that for every $p \in Per_H \cap B_\nu(x)$ we have that $S_p \in W^s_\eps(x)$ where $0<2\eps< \alpha$ and $\alpha>0$ is the constant of expansivity.

\begin{lema}\label{bordeenestable} For every maximal
chain $[p]$,  $\partial \left(\bigcup_{q\in [p]} \umb{q}\right)
\cap B_\nu (x)\subset W^s_\eps(x)$ verifies.
\end{lema}

\dem{} Because of Lemma \ref{clasesdisjuntas} we know that $[p]$
accumulates in $x$. Let $q_n \in [p]$ such that $q_n \to x$. Take
a point $y \in \partial \left(\bigcup_{q\in [p]} \umb{q}\right)$. Then, a sequence $z_n \in
\umb{q}_n$ exists such that $z_n \to y$. Without loss of
generality we can suppose $z_n \in S_{q_n}$.

Remark \ref{observacionestablebola} ensures the existence of $\eps>0$
such that $S_{q_n} \subset W^s_{\eps}(z'_n)$. So we have that for
all $m\geq 0$

$$\dist(f^{m}(y),f^{m}(x))= \lim_{n\to \infty} \dist(f^m (z'_n),f^m (p_n))
\leq \eps$$

\noindent so $y \in W^s_{\eps}(x)$. Then, $\partial \left(\bigcup_{q\in [p]} \umb{q}\right)\subset
W^s_\eps(x)$.

\lqqd

\begin{lema}\label{cub} Suppose $\overline{Per_H}=M$ and let $x\in M$ be a singularity.
  Then, for all $p\in B_\delta (x)\cap Per_H^{n-1}$, there exists
  a neighborhood $V$ of $S_p$ such that
  $Per_H \cap V \cap B_\delta (x)\subset [p]$.
\end{lema}

\dem By contradiction, let us suppose that there exists $y \in S_p\cap B_{\delta}(x)$
satisfying that $q_n\to
  y$ with $q_n\in [q]\neq [p]$ (remember that because of Theorem \ref{teo1},
  near $S_p$ we have local product structure
  so every periodic point near $y$ must have the same index as $p$).
  Then $y\in \partial S_{[q]}$ (because it belongs both to $int(S_{[p]})$ and $S_{[q]}$, and the interiors of $S_{[p]}$ and $S_{[q]}$ have empty intersection).
  So, by Lemma \ref{bordeenestable},
  $y\in W^s_{\eps}(x)$. Therefore  $x\in W^s(p)$ because $y \in S_p \en W^s_{\eps}(p)$. But, since
  $p\in Per_H$ we contradict the fact that $x$ is singular, since Theorem \ref{teo1}
  gives us local product structure in a neighborhood of $x$ by iteration of the local product structure in $p$.

\lqqd

\begin{lema}\label{noespinas}
  If $\overline{Per_H}=M$ and let $x\in M$ be a singularity. Then
  $\interior(S_{[p]})\cap B_\delta (x)=\bigcup_{q\in [p]} \umb{q}\cap B_\delta(x)$.
\end{lema}

\dem

  The inclusion $\bigcup_{q\in [p]} \umb{q}\cap B_\delta(x)\en \interior(S_{[p]})\cap B_\delta (x)$ is immediate because if $q \geq r$ then $\umb{r} \en int(\umb{q})$.

  To obtain the other inclusion we proceed by contradiction
  supposing there exists a point $y \in B_\delta (x)$ in the
  interior of $S_{[p]}$ but such that $y\notin \umb{q}$ for all $q
  \in [p]$.

  Then, there exists $y_n\in S_{q_n}$ such that $y_n \to
  y$ (this implies in particular that $y \in W^s_\eps(x)$ because of Lemma \ref{semicontinuidad}) where $q_n \in [p]$ satisfies $q_n \to x$.

  Using Lemma \ref{cub} and the fact that $\overline{Per_H}=M$ we know that there exist points $r_n \in [p]$  arbitrarily close to $y_n$. We can suppose $r_n \to
  y$ and that this points are not bounded in the ordering in $[p]$. On the other hand, we consider $U_{r_n}=cc_{r_n}(B_\nu(x)
  \cap W^u(r_n)) \en W^u_\eps(r_n)$ (see Remark \ref{observacionestablebola}) which is separated by $S_{r_n}$ in two different connected components (see corollary \ref{sepcoduno}).

  Pick  $\gamma>0$ and choose $z_n \in
  \partial B_\gamma(y)\cap U_{r_n}$ such that $z_n \notin \umb{r_n}$.
  We can suppose that $z_n \to z \in \partial B_\gamma (y)$ and
  using the semicontinuous variation of local stable and unstable
  sets (Lemma \ref{semicontinuidad}) we obtain that $z \in
  W^u_{\eps}(y)$.

  We shall prove that $z \notin S_{[p]}$ and since $\gamma$ was arbitrary this will imply that $y \in \partial S_{[p]}$ which contradicts the fact that $y \in int(S_{[p]})$.

  We know that $z \notin \umb{q}$ for all $q \in [p]$, so, if $z \in S_{[p]}$ it should be
  accumulated by points in $S_{q_n}$ and therefore verify $z \in
  W^s_{\eps}(x)$. Then, $z \neq y$, $z\in W^u_{\eps}(y)$ and $y,z\in W^s_{\eps}(x)$
  which contradicts expansivity (remember we chose $\eps$ so that $2\eps<\alpha$).

\lqqd

\begin{obs}\label{bordesepara}
  Clearly $x\in S_{[p]}$ and $x\notin \umb{q}$ for all
  $q\in [p]$. Since $S_{[p]}$ is a closed set with non empty interior and $x\in\partial S_{[p]}$ we have that its complement in $B_\nu(x)$ which is open
  is also non empty. This implies that $\partial S_{[p]}$
  separates $B_\nu (x)$.
\end{obs}

The next lemma shows how the stable sets of periodic points converge uniformly towards $\partial S_{[p]}$.

\begin{lema}\label{lemabordepar}
  Suppose $\overline{Per_H}=M$ and let $z \in \partial S_{[p]} \cap
  B_{\delta}(x)$ and $\rho >0$. Then, there exists $V$ a neighborhood of $z$ such that
  if $q \in [p] \cap V$ then $S_{[p]}\cap
  \overline{B_{\delta}(x)} \en \umb{q} \cup B_{\rho}(\partial
  S_{[p]})$.
\end{lema}

\dem{} Given $\rho>0$, the set $K =\left(S_{[p]} \setminus
B_{\rho}(\partial S_{[p]})\right) \cap \overline{B_{\delta}(x)}$
is a compact set contained in $\interior(S_{[p]}\cap B_\delta
(x))$ so, using Lemma \ref{noespinas}, $\{\interior(\umb{q})\}_{q \in
[p]}$ is an open cover of $K$ so $r\in [p]$ exists such that
$K\subset \umb{r}$. Let $V$ be a neighborhood of $z$ disjoint from
$\umb{r}$. Then, for every $q\in [p]\cap V$ we have that $q\geq
r$. Then, $K=\left(S_{[p]} \setminus B_{\rho}(\partial
S_{[p]})\right) \cap \overline{B_{\delta}(x)}\subset \umb{q}$ and
therefore $S_{[p]}\cap
  \overline{B_{\delta}(x)} \en \umb{q} \cup B_{\rho}(\partial
  S_{[p]})$.

\lqqd


 The following Lemma represents the key step for proving
the uniformity of the local product structure because it allows us
to ensure that the stable and unstable sets intersect in a
neighborhood of a singularity. This gives uniformity and is also
important  to give structure to $\partial S_{[p]}$ and discard
singularities.

\begin{lema}\label{propinterseccion} Suppose $\overline{Per_H}=M$.
For all
 $z \in \partial S_{[p]}\cap B_\delta (x)$ and for all
 $\eps>0$ there exists $V$ neighborhood of $z$ such that if $q, r \in V \cap [p]$ then $U_q$
intersects $S_r$ and $\partial S_{[p]}$ in $B_\eps(z)\cap
S_{[p]}$.
\end{lema}

\dem{} Let $V$ be a neighborhood of $z$ such that $z\in V\subset
B_{\delta}(x)$. Corollary \ref{sepcoduno} allows us to associate
to each $q\in [p]\cap V$ two points $y_1^q, y_2^q\in U_q\cap
\partial B_\delta(x)$ such that $y_1^q\in
\umb{q}$ and $y_2^q \notin \umb{q}$. Lemma \ref{angulos} gives us
$\rho>0$ such that (maybe taking $V$ smaller) for $i=1,2$ and
$q\in [p]\cap V$

\begin{equation}\label{propinterseccionEc1}
\dist(y_i^q, \partial S_{[p]})>\rho
\end{equation}

At the same time, by Lemma \ref{lemabordepar} we can suppose that
for every $q\in[p]\cap V$,
\begin{equation}\label{propinterseccionEc2}
  S_{[p]} \cap B_{\delta}(x) \en
  \umb{q} \cup B_{\rho/2}(\partial S_{[p]})
\end{equation}

Then, since $y_2^q\notin \umb{q}$ and $y_2^q\notin B_\rho(\partial
S_{[p]})$, we have that $y_2^q\notin S_{[p]}$. Remark
\ref{bordesepara} together with the fact that $U_q$ is connected
implies $U_q$ intersects $\partial S_{[p]}$.

Let us take $r \in V\cap [p]$ such that $q\leq r$, that is to say
$\umb{q}\subset \umb{r}$. Consider $y_1^r$ and $y_2^r$ associated
to $r$ in the same way we did with $q$. Then $y_1^q\in
\umb{q}\subset \umb{r}$. On the other hand, $y_1^r\in S_{[p]}\cap
\partial B_\delta(x)$ and by (\ref{propinterseccionEc2}) $y_1^r\in
\umb{q} \cup B_{\rho/2}(\partial S_{[p]})$. Because of
(\ref{propinterseccionEc1}) we have $y_1^r\notin
B_{\rho/2}(\partial S_{[p]})$ and so $y_1^r\in \umb{q}$. Then
$y_1^r,y_1^q\in \umb{q}\subset \umb{r}$.

Now, since $y_2^r\notin \umb{r}$ and $\umb{q}\subset \umb{r}$ we
have that $y_2^r\notin \umb{q}$. Previously we said that
$y_2^q\notin S_{[p]}$, so, applying (\ref{propinterseccionEc2})
(to $r$ instead of $q$) we have that $y_2^q\notin\umb{r}$. This
implies $y^q_2,y^r_2\notin \umb{r}\supset\umb{q}$.

Finally, since $U_q\supset\{y_1^q,y_2^q\}$ and
$U_r\supset\{y_1^r,y_2^r\}$ are connected, and $S_q$ and $S_r$
separate the ball $B_\nu (x)$ we deduce that $S_q\cap U_r$ and
$U_q\cap S_r$ are not empty as wanted.

Given $\eps >0$, expansivity and semicontinuous variation of local
stable and unstable sets allow us to prove that by means of
considering $V$ small enough we can ensure that the intersections
lie in $B_\eps(z)$.

\lqqd

To prove Theorem \ref{teo2} we shall also make use of some
properties of the frontier of the sets $S_{[p]}$.

\begin{prop}\label{lem3dim} If $\overline{Per_H}=M$, $\partial S_{[p]} \cap
B_\delta(x)$ is a topological manifold of dimension $n-1$.
\end{prop}

\dem{} Let $z\in \partial S_{[p]} \cap B_\delta(x)$. We choose
$\eps>0$ such that $B_{\eps}(z) \en B_{\delta}(x)$ and let $V$ a
neighborhood of $z$ satisfying that if $q, r \in V \cap [p]$ then
$U_q \cap S_r \cap B_{\eps}(z) \neq \emptyset$ as given in  Lemma \ref{propinterseccion}. Also, for every $q \in V \cap [p]$ we can have
$U_q \cap
\partial S_{[p]} \cap B_{\eps}(z)\neq \emptyset$ again by Lemma \ref{propinterseccion}.

Pick $q \in V \cap [p]$ and define $h_q\colon S_q \cap
\overline{V} \to
\partial S_{[p]} \cap \overline{B_{\eps}(z)}$ given by

$$h_q(y)=
\lim_{q_n \to y} U_{q_n} \cap \partial S_{[p]}$$

\noindent which is well defined thanks to expansivity and
semicontinuous variation of local stable and unstable sets (Lemma
\ref{semicontinuidad}) together with the fact that $\partial
S_{[p]} \en W^s_{\gamma}(z)$ because of Lemma \ref{lemabordepar}.
The fact that there is a sequence $q_n \in [p] \to y$ is a consequence of Lemma \ref{cub} and the fact that $\overline{Per_H} =M$.

The same argument implies that $h_q$ is continuous and injective.
Moreover, since the domain is compact, $h_q$ is a homeomorphism over
its image.

Again, by Lemma \ref{propinterseccion} , we can take $V'$ and $\eps'>0$ such that $\overline {B_{\eps'}(z)} \en V$ and that for $q,r \in [p] \cap V'$, $U_q \cap S_r \cap
B_{\eps'}(z) \neq \emptyset$. Analogously, we have that for
every $q \in V' \cap [p]$,  $U_q \cap \partial S_{[p]} \cap
\overline{B_{\eps'}(z)}\neq \emptyset$ verifies.

If we fix $q \in V' \cap [p]$ (Lemmas \ref{cub} and \ref{lemabordepar} and ensures the existence of such $q$) we will be able to prove that for all $w \in \partial S_{[p]} \cap V'$ exists $y \in S_q \cap V$ such that $h_q(y)=w$. This holds since for
every $w \in \partial S_{[p]} \cap V'$ we can find $\{q_n\} \en
[p] \cap V'$ such that $q_n \to w$ and so that $\emptyset \neq
U_{q_n} \cap S_q \cap B_{\eps'}(z) \en V \cap S_q$. In particular,
every point in $\partial S_{[p]} \cap V'$  has a preimage of the
map $h_q$ in  $S_q \cap \overline{ B_{\eps'}(z)} \en S_q \cap V$.

Since $h_q$ is a homeomorphism over its image, $\partial
S_{[p]} \cap V'$ is homeomorphic to its preimage which is an open subset of $S_q \cap V'$ and the proposition is proved (remember $S_q \cap V$ is homeomorphic to an open
set of $\R^{n-1}$).
\lqqd


\subsection{Constant splitting}\label{secciondimconst}

\demos{of Theorem \ref{teo2}}

By contradiction, we suppose that
$\emptyset\neq\overline{Per_H^{n-1}}\neq M$ and consider a
singularity  $x \in
\partial \overline{Per_H^{n-1}}$.
We consider $\nu$ and $\delta$ as in Lemma \ref{clasesdisjuntas},
for which we know there is a finite set of maximal chains of the
partial order in $X_\delta$. Let $[p]$ be a maximal chain
accumulating in $x$.

\begin{lema}\label{densoenparaguas}
There exists $\delta>0$ such that $Per_H \cap S_{[p]}\cap
B_\delta(x) \en Per_H^{n-1}$.
\end{lema}
\dem Suppose, by contradiction, that there exist $p_n,q_n\to x$
where $q_n\in S_{[p]}\cap Per_H\setminus Per_H^{n-1}$ and $p_n\in
[p]$.

We know that $p_n,q_n\notin \partial S_{[p]}$ because it would contradict the fact that $x$ is a singularity.

Since $q_n \notin Per_H^{n-1}$, $U_{q_n}$ is a connected topological manifold
(and therefore arcconnected) of dimension at least two.
Consequently, if we remove a point from $U_{q_n}$ it remains
arcconnected.

Clearly, for every $p_n$ there exists $q_m\notin
\umb{p_n}$. Remember that $\partial S_{[p]}$ and $S_{p_n}$
separates the ball $B_\nu(x)$.

We shall prove that $U_{q_m}\subset S_{[p]}\setminus \umb{p_n}$.
Otherwise, $y\in U_{q_m}\setminus S_{[p]}$ would exist. Since
$\partial S_{[p]}$ separates the ball $B_\nu(x)$ we know that
every curve contained in $U_{q_m}$ joining $q_m$ to $y$ must
intersect $\partial S_{[p]}$. Expansivity implies that $U_{q_m}$
intersects $\partial S_{[p]}$ in at most one point. Then, since
two curves in $U_{q_m}$ connecting $q_m$ to $y$ and coinciding
only in the extremes exist (because of the dimension of $U_{q_m}$)
they should intersect $\partial S_{[p]}$ in two different points
reaching a contradiction. We proceed analogously if we consider
$y\in \umb{p_n}$.

Finally, the fact that for every $n_0$ there exist $m,n\geq n_0$ such
that $U_{q_m}\subset S_{[p]}\setminus \umb{p_n}$ contradicts
expansivity (see Lemma \ref{angulos}).

\lqqd

Let $\cmax$ be the finite set of maximal chains in $B_\delta (x)$
and let
\[
  S=\bigcup_{[p] \in \cmax} S_{[p]}
\]
Since every $S_{[p]}$ is closed in $B_\nu(x)$ and $\cmax$ is
finite, we have that $S$ is closed. Lemma \ref{densoenparaguas} and the fact that $\overline{Per_H}=M$ implies

\begin{equation}\label{paraguero}
  B_{\delta}(x)\cap S = B_{\delta}(x) \cap \overline{Per_H^{n-1}}
\end{equation}

Since $x \in
\partial \overline{Per_H^{n-1}}$
we know that $S$ can not be a neighborhood of $x$. We shall see
how this fact represents a contradiction.

In order to do that, we shall make use of Proposition
\ref{lem3dim} and the following lemma.

\begin{lema}\label{lem2dim} For all
$p \in Per_H^{n-1}\cap B_\delta(x)$ exists $A_{[p]}\subset\partial
S_{[p]}$ such that $A_{[p]}$ is an open and dense subset relative
to $\partial S_{[p]}\cap B_\delta(x)$ and $A_{[p]}$ is in the
interior of $S$.
\end{lema}

Proposition \ref{lem3dim} ensures that $\partial S_{[p]}\cap
B_\delta(x)$ is a topological manifold of dimension $n-1$. Then,
 Lemma \ref{lem2dim} and a result in \cite{HW} stating that a closed set with empty
interior in a topological manifold has dimension smaller than the
manifold (chapter IV, section 4), imply that for every $[p]$, $\dimtop(\partial
S_{[p]}\setminus A_{[p]})\leq n-2$. Moreover, since  $\partial S
\en \bigcup \partial S_{[p]}$

$$\partial S\subset \bigcup_{[p]\in
\cmax} \partial S_{[p]}\setminus
A_{[p]}$$

And, since the union of a finite set of closed spaces has the
dimension of the largest one (see \cite{HW} chapter III,
section 3) we know that  $\dimtop(\partial S)\leq n-2$. So,
$\partial S$ can not separate $B_\delta(x)$ because it should have
dimension at least $n-1$ (see \cite{HW} chapter IV, section 5).
This leads us to a contradiction.

\demos{of Lemma \ref{lem2dim}}

Let $\eps>0$ and $z\in\partial S_{[p]}$. By Lemma
\ref{propinterseccion}, there exist $q \in [p]$ such that
$\{a\}=U_q \cap \partial S_{[p]}$ is in $B_{\eps}(z)$. Theorem
\ref{teo1} implies that $q$ has a neighborhood with local product
structure, by iterating this neighborhood to the past, we obtain
local product structure over a neighborhood of $a$, so, $a$ must belong to $int(S)= int(\overline{Per_H^{n-1}}\cap B_{\delta}(x))$  and the lemma is proved.

\lqqd

\subsection{Uniform local product structure}\label{seccionsingularidades}

We shall prove Theorem \ref{teo3} in this section. By Theorem
\ref{teo1} we know that there is an open and dense set whose
points admit a local product structure. And by Theorem \ref{teo2}
we conclude, since $Per_H^{n-1}\neq\emptyset$, that
$Per_H=Per_H^{n-1}$.

Let $\sing$ be the set of singularities of $f$, that is to say,
the points which do not admit any local product structure. To
prove Theorem \ref{teo3} we must prove that $\sing$ is an empty
set.

With the results proved in \ref{seccionlemas} we obtain the following
consequence which allows us to study the set of singularities in
codimension one case. The next proposition gives a sort of local product
structure in the sets $S_{[p]}$ which will be defined properly in
this statement.

\begin{prop}\label{paragualocal} Let $x \in \sing$. Then, for every $z \in
\partial S_{[p]} \cap B_{\delta}(x)$ there exists $h:I \times I^{n-1}
\to S_{[p]}$ ($I = [0,1]$) homeomorphism over its image, where
$h(\{a\} \times I^n)$ is contained in a local stable set, $h(I
\times \{b\})$ is contained in a local unstable set and the image
of $h$ is a neighborhood of $z$ relative to $S_{[p]}$.
\end{prop}

\dem{} By Lemma \ref{propinterseccion} there exists $V \en
B_{\delta}(x)$ neighborhood of $z$ in $M$ such that if $q, r \in
[p] \cap V$ then $S_q \cap U_r \neq \emptyset$ and $U_q \cap \partial S_{[p]}\neq \emptyset$.

Let $D_z \en \partial S_{[p]} \cap V$ homeomorphic to $I^{n-1}$
(see Proposition \ref{lem3dim}) such that $z$ belongs to the
interior of $D_z$ relative to $\partial S_{[p]}$.

Let $V' \en V$ neighborhood of $z$ such that if $q \in [p] \cap
V'$ then $S_q \cap U_r \cap V \neq \emptyset$ and $U_q \cap D_z \neq \emptyset$.

Let $q \in V' \cap [p]$ and we define $h: U_q \cap \overline V \cap S_{[p]} \times D_z \to
S_{[p]}$ in such a way that $h(y,w) = W^s_{\eps}(y) \cap
W^u_{\eps}(w)$ is verified. By the choice of $V$, approximating
with topologically hyperbolic periodic points and making use of
expansivity and semicontinuous variation of local stable and
unstable sets (Lemma \ref{semicontinuidad}) we can ensure that the
map $h$ is well defined, continuous and injective and, since the
domain is compact, a homeomorphism over its image.

We are now interested in proving that the image contains $V'\cap
S_{[p]}$ and it is enough to show that it contains $[p] \cap V'$,
since $Per_H^{n-1}$ is dense in $V''$.  This holds due to the choice of $V'$.

Let $U = h^{-1}(V' \cap S_{[p]})$ which is open because $h$ is a
homeomorphism over its image. Since $z \in V' \cap S_{[p]}$ a
relative open set of the image of $h$ and $S_{[p]}$, $h^{-1}(z)$
is in the interior of $U$. Since $U_q \cap \overline V \cap
S_{[p]} \times D_z$ is locally connected in $h^{-1}(z)$ we can
find in $U$ a set homeomorphic to $I \times I^{n-1}$ neighborhood
of $h^{-1}(z)$ whose image will be a relative neighborhood of $z$
in $S_{[p]}$.

The other properties of $h$ claimed in the statement of the proposition are immediate consequences of the definition of $h$.

\lqqd

\begin{lema}\label{sing2}
  If $\overline{Per_H^{n-1}}=M$ then $\sing$ is a finite set.
\end{lema}

\dem{}

Since the set of points with local product structure is open and
invariant we know that $\sing$ is compact and invariant.
Therefore, $f\colon \sing\to\sing$ is an expansive homeomorphism.

 We shall prove that there exist a neighborhood
of $x \in \sing$ satisfying that every singularity in that neighborhood
belongs to the local stable set of $x$. This is a consequence of
the existence of $\delta>0$ small enough (given by Lemma
\ref{clasesdisjuntas}) such that (since $Per_H^{n-1}$ is dense) we
have that  $B_{\delta}(x) \subset \bigcup_{i=1}^{k} S_{[p_i]}$.
Proposition \ref{paragualocal} implies that in the interior of
$S_{[p_i]}$ there is a local product structure (maybe by
considering $\delta$ smaller) so singularities must lie in
$\bigcup_{i=1}^{k}
\partial S_{[p_i]}$. Lemma \ref{bordeenestable} now implies that singularities of
$B_\delta (x)$ belong to the local stable set of $x$.

Expansivity implies that Lyapunov stable points are asymptotically
stable. Otherwise, points $y,w$ such that $\dist(f^n(y), f^n(w))
\leq \eps \leq \alpha$ ($\alpha$ expansivity constant) and such
that a subsequence $n_j \to +\infty$ with $\dist(f^{n_j}(y),
f^{n_j}(w)) \geq \delta$ exist. Taking limit points we contradict
expansivity.

Since $\sing$ is compact and every point is asymptotically stable
for $f$, we conclude that $\sing$ must be finite.

\lqqd

In the following Lemma we will show that there are no isolated singularities if $\dim(M)\geq 3$. Observe that in surfaces, pseudoAnosov maps have this kind of singular points. The key fact here is how the semilocal product structures given by Proposition \ref{lem3dim} are glued around the singularity. The idea is that if $S_{[p]}$ and $S_{[q]}$ have semilocal product structure, then $S_{[p]}\cap S_{[q]}\setminus \{x\}$ is a connected component of $W^s_{loc}(x)\setminus \{x\}$. If $\dim(M)\geq 3$ then the set $W^s_{loc}(x)\setminus \{x\}$ is connected and therefore there is no place for a third semilocal product structure. This will let us prove that $x$ has a local product structure.

\begin{lema}\label{sing1}
    If $\dim (M)\geq 3$ and $\overline{Per_H^{n-1}}=M$ then, no
    isolated singularities exist.
\end{lema}

\dem
  By contradiction, suppose $x\in M$ is an isolated singularity.
  Let $\nu,\delta>0$ be as in Proposition \ref{lem3dim} and such that
  $B_\nu(x)\cap\sing=\{x\}$.
  Fix $[p]$ a maximal chain accumulating in $x$ and let $T=cc_x(\partial S_{[p]}\cap B_\delta(x))$. We know by Proposition \ref{lem3dim}  that $T$ is a topological manifold that is closed in $B_\delta(x)$.

\par Let $z\in T\setminus\{x\}$ and $[q]\neq [p]$ such that $z\in\partial S_{[q]}$. Define $T'=cc_z(\partial S_{[q]}\cap B_\delta (x))$. Let $F=T'\setminus\{x\}\cap T\setminus\{x\}$, which is a non empty closed set in both $T\setminus\{x\}$ and $T'\setminus\{x\}$. Since for all $w\in F$ there is a local product structure, $F$ is open in both $T\setminus\{x\}$ and $T'\setminus\{x\}$. Thus $F=T\setminus \{x\}=T'\setminus\{x\}$ because $T\setminus \{x\}$ and $T'\setminus\{x\}$ are connected sets ($\dim(M)\geq 3$). Then, since $T'$ is closed, $x \in T'$ which implies that $T'=cc_x(\partial S_{[q]} \cap B_\delta(x))$.

  Since $Per_H^{n-1}$ is dense in $S_{[p]}$
  we can apply Proposition \ref{paragualocal} to $x$. Let
  $$
    h_p:[0,1) \times (-1,1)^{n-1} \to R_p\subset B_\delta(x)
  $$
  \noindent be a homeomorphism such
  that $R_p$ is a neighborhood of $x$ relative
  to $S_{[p]}$ and $h_p(0)=x$. Let $F_p=T\cap R_p=h(\{0\} \times (-1,1)^{n-1})$.

  Now, from Proposition \ref{paragualocal} we can consider
  $h_q:(-1,0] \times (-1,1)^{n-1} \to R_q\subset B_\delta(x)$
  a homeomorphism satisfying that $R_q$ is a neighborhood of $z$
  relative to $S_{[q]}$ and $h_q(0)=x$. Analogously we define
  $F_q=\partial S_{[q]}\cap R_q=h_q(\{0\} \times (-1,1)^{n-1})\subset F$.
  From the previous, we can suppose $F_q\subset F_p$.
  Let
  $\pi_2 \colon \R \times \R^{n-1}\to \R^{n-1}$ the canonical projection over the second
  coordinate.
  Furthermore, if we restrict $h_p$ to the set $[0,1)\times \pi_2(h_q^{-1}(F_q))$
  we can suppose $F_p=F_q$.

  Let $h\colon (-1,1)\times F_p\to
  B_\delta(x)$ given by

  $$
    h(t,y)=
    \left\{
      \begin{array}{ll}
        h_p(t,\pi_2(h_p^{-1}(y)))& \hbox{if $t\geq 0$}
        \\
        h_q(t,\pi_2(h_q^{-1}(y)))& \hbox{if $t\leq 0$}
        \\

      \end{array}
    \right.
  $$

  Clearly $h(0,y)=y$ so $h$ is continuous. Again, using the Invariance of
  Domain Theorem, this allows us to prove that $h$ gives a
  local product structure around $x$.
  This contradicts the fact that $x$ is a singularity.
\lqqd

\demos{of Theorem \ref{teo3}}

\medskip

Once we have discarded singularities it is very simple to prove
there is a uniform local product structure. Otherwise, there would
exist points $x_n$ not admitting local product structure in balls
of radius greater than $1/n$. Taking a limit point we could find a
singularity, a contradiction.

Uniform local product structure implies the pseudo orbit tracing
property from the results of \cite{R} which ensure the
existence of a hyperbolic metric in the coordinates given by the
local product structure (see \cite{V2}).

\lqqd


\section{Appendix}\label{secciontoro}

To conclude, we prove the following proposition and then sketch
the proof that $M$ is $\T^n$.

\begin{prop}\label{estprodloc} Let $M$ be a $n-$dimensional manifold ($n\geq 3$) and
$f:M \to M$ an expansive homeomorphism such that $Per_H$ is dense in
$M$ and $Per_H^1\neq \emptyset$ or $Per_H^{n-1}\neq \emptyset$.
Then, $M$ admits a codimension one foliation with leaves
homeomorphic to $\R^{n-1}$.
\end{prop}

\dem{} The uniform local product structure obtained in Theorem \ref{teo3}
shows the existence of the foliation.

Let us suppose that $Per_H^{n-1}\neq \emptyset$, then, the leaves of
the foliation are the stable sets of the points. Let $x\in M$, we
shall prove that $W^{s}(x)$ is homeomorphic to $\R^{n-1}$. To see this,
is enough to see that

\[ W^s(x) = \bigcup_{n\geq 0} f^{-n} ( S_{\eps}(f^{n}(x)) \]

Where $S_\eps(z)$ is a disc of uniform size in $W^s_{\eps}(z)$
(which exist because of the uniform local product structure). So,
$W^s(x)$ may be written (maybe by taking some subsequence $n_j \to
\infty$ so that $f^{-n_j}(S_{\eps}(f^{n_j}(x))) \subset
f^{-n_{j+1}} (S_{\eps}(f^{-n_{j+1}}(x)))$) as an increasing union
of $n-1$ dimensional discs, which implies the thesis.

\lqqd

Once we know the leaves are homeomorphic to $\R^{n-1}$ classical
arguments allow us to prove that $M$ is $\T^n$. As we said, we
shall sketch some steps of the proof for the sake of completeness.
The ideas are  based on \cite{V3} and
\cite{F} section $5$.

The first thing it should be proved is that the universal covering
space of $M$ ($\overline{M}$) equals $\R^n$.

To prove that $\overline{M}=\R^{n}$ it suffices to prove that
given two points $\overline{x}, \overline{y} \in \overline{M}$
then, the lifts of their stable and unstable manifolds (which are
respectively proper copies of $\R^{n-1}$ and $\R$) intersect at a
single point.

To see that the intersection has at most one point, we can
see that if the manifolds intersect at more than one point then we
can obtain a closed loop transversal to the codimension one
foliation, thus, bounding a disc (since we are in the universal
covering, the loop is nullhomotopic). By using Solodov's methods
(see \cite{So} Lemma 5) we see that the disc may be chosen to
be in general position so that we obtain a foliation of the disc
$\D^2$, transversal to the frontier and such that its
singularities are nondegenerate and have no saddle connections
(this is the only step where differentiability is used in
\cite{F}). Now, using Haefliger arguments (see
\cite{V3} Lemma 2.11 or \cite{F} Lemma 5.1) we
conclude there is a leaf of the codimension one foliation with non
trivial holonomy, hence, the leaf is not simply connected, a
contradiction.

Finally, proving that the foliations intersect is a straightforward
adaptation of the arguments of \cite{F} Lemma 5.2 after it is
known that the the leaves of the codimension one foliation are
dense (which follows from the fact that periodic points are dense
and the uniform local product structure).

Once this is obtained, it is not difficult to prove that $\pi_1(M)$ is free abelian by studding the action of $\pi_1(M)$ over $\R$ as it permutes without fixed points the leaves of the codimension one foliation (see \cite{HeHi} Chapter VIII, section 3, remember that
the leaves of the foliation are dense).

Now one can follow the proof in \cite{F},
by reading the proofs of Proposition (6.2), Theorem (4.2) and
Theorem (3.6) in that order (remember that expansive
homeomorphisms with local product structure have hyperbolic
canonical coordinates, \cite{R}).

One can take a shortcut in dimensions $\geq 5$ thanks to a result of \cite{HiWa}.  A space with free abelian fundamental group and which is covered by $\R^n$ is an
Elienberg-McLane space of the same type of a torus, hence
homotopically equivalent to one (see \cite{Hat}, Theorem
1.B.8.). From \cite{HiWa} we deduce that if $n$, the dimension
of $M$, satisfies $n\geq 5$ then $M$ is homeomorphic to $\T^n$.

This proves that $M=\T^n$.

\lqqd


\end{document}